\documentclass[12pt]{article}
\usepackage{amsmath,amssymb,amsbsy,amsfonts,amsthm,latexsym,amsopn,amstext,
amsopn,amstext,amsxtra,euscript,amscd}
\begin{document}
\newfont{\teneufm}{eufm10}
\newfont{\seveneufm}{eufm7}
\newfont{\fiveeufm}{eufm5}
%
%
\newfam\eufmfam
                                 \textfont\eufmfam=\teneufm
\scriptfont\eufmfam=\seveneufm
                                 \scriptscriptfont\eufmfam=\fiveeufm

%
%
\def\frak#1{{\fam\eufmfam\relax#1}}
%


\def\rad{{\rm rad}}

\def\bbbr{{\rm I\!R}} 
\def\bbbc{{\rm I\!C}} 
\def\bbbm{{\rm I\!M}}
\def\bbbn{{\rm I\!N}} 
\def\bbbf{{\rm I\!F}}
\def\bbbh{{\rm I\!H}}
\def\bbbk{{\rm I\!K}}
\def\bbbl{{\rm I\!L}}
\def\bbbp{{\rm I\!P}}
\def\bbbq{{\rm I\!Q}}
\newcommand{\lcm}{{\rm lcm}}
\def\bbbone{{\mathchoice {\rm 1\mskip-4mu l} {\rm 1\mskip-4mu l}
{\rm 1\mskip-4.5mu l} {\rm 1\mskip-5mu l}}}
\def\bbbc{{\mathchoice {\setbox0=\hbox{$\displaystyle\rm C$}\hbox{\hbox
to0pt{\kern0.4\wd0\vrule height0.9\ht0\hss}\box0}}
{\setbox0=\hbox{$\textstyle\rm C$}\hbox{\hbox
to0pt{\kern0.4\wd0\vrule height0.9\ht0\hss}\box0}}
{\setbox0=\hbox{$\scriptstyle\rm C$}\hbox{\hbox
to0pt{\kern0.4\wd0\vrule height0.9\ht0\hss}\box0}}
{\setbox0=\hbox{$\scriptscriptstyle\rm C$}\hbox{\hbox
to0pt{\kern0.4\wd0\vrule height0.9\ht0\hss}\box0}}}}
\def\bbbq{{\mathchoice {\setbox0=\hbox{$\displaystyle\rm
Q$}\hbox{\raise
0.15\ht0\hbox to0pt{\kern0.4\wd0\vrule height0.8\ht0\hss}\box0}}
{\setbox0=\hbox{$\textstyle\rm Q$}\hbox{\raise
0.15\ht0\hbox to0pt{\kern0.4\wd0\vrule height0.8\ht0\hss}\box0}}
{\setbox0=\hbox{$\scriptstyle\rm Q$}\hbox{\raise
0.15\ht0\hbox to0pt{\kern0.4\wd0\vrule height0.7\ht0\hss}\box0}}
{\setbox0=\hbox{$\scriptscriptstyle\rm Q$}\hbox{\raise
0.15\ht0\hbox to0pt{\kern0.4\wd0\vrule height0.7\ht0\hss}\box0}}}}
\def\bbbt{{\mathchoice {\setbox0=\hbox{$\displaystyle\rm
T$}\hbox{\hbox to0pt{\kern0.3\wd0\vrule height0.9\ht0\hss}\box0}}
{\setbox0=\hbox{$\textstyle\rm T$}\hbox{\hbox
to0pt{\kern0.3\wd0\vrule height0.9\ht0\hss}\box0}}
{\setbox0=\hbox{$\scriptstyle\rm T$}\hbox{\hbox
to0pt{\kern0.3\wd0\vrule height0.9\ht0\hss}\box0}}
{\setbox0=\hbox{$\scriptscriptstyle\rm T$}\hbox{\hbox
to0pt{\kern0.3\wd0\vrule height0.9\ht0\hss}\box0}}}}
\def\bbbs{{\mathchoice
{\setbox0=\hbox{$\displaystyle     \rm S$}\hbox{\raise0.5\ht0\hbox
to0pt{\kern0.35\wd0\vrule height0.45\ht0\hss}\hbox
to0pt{\kern0.55\wd0\vrule height0.5\ht0\hss}\box0}}
{\setbox0=\hbox{$\textstyle        \rm S$}\hbox{\raise0.5\ht0\hbox
to0pt{\kern0.35\wd0\vrule height0.45\ht0\hss}\hbox
to0pt{\kern0.55\wd0\vrule height0.5\ht0\hss}\box0}}
{\setbox0=\hbox{$\scriptstyle      \rm S$}\hbox{\raise0.5\ht0\hbox
to0pt{\kern0.35\wd0\vrule height0.45\ht0\hss}\raise0.05\ht0\hbox
to0pt{\kern0.5\wd0\vrule height0.45\ht0\hss}\box0}}
{\setbox0=\hbox{$\scriptscriptstyle\rm S$}\hbox{\raise0.5\ht0\hbox
to0pt{\kern0.4\wd0\vrule height0.45\ht0\hss}\raise0.05\ht0\hbox
to0pt{\kern0.55\wd0\vrule height0.45\ht0\hss}\box0}}}}
\def\bbbz{{\mathchoice {\hbox{$\sf\textstyle Z\kern-0.4em Z$}}
{\hbox{$\sf\textstyle Z\kern-0.4em Z$}}
{\hbox{$\sf\scriptstyle Z\kern-0.3em Z$}}
{\hbox{$\sf\scriptscriptstyle Z\kern-0.2em Z$}}}}
\def\ts{\thinspace}

\newtheorem{theorem}{Theorem}
\newtheorem{lemma}[theorem]{Lemma}
\newtheorem{claim}[theorem]{Claim}
\newtheorem{cor}[theorem]{Corollary}
\newtheorem{prop}[theorem]{Proposition}
\newtheorem{definition}{Definition}
\newtheorem{question}[theorem]{Open Question}

\def\squareforqed{\hbox{\rlap{$\sqcap$}$\sqcup$}}
\def\qed{\ifmmode\squareforqed\else{\unskip\nobreak\hfil
\penalty50\hskip1em\null\nobreak\hfil\squareforqed
\parfillskip=0pt\finalhyphendemerits=0\endgraf}\fi}

\def\cA{{\mathcal A}}
\def\cB{{\mathcal B}}
\def\cC{{\mathcal C}}
\def\cD{{\mathcal D}}
\def\cE{{\mathcal E}}
\def\cF{{\mathcal F}}
\def\cG{{\mathcal G}}
\def\cH{{\mathcal H}}
\def\cI{{\mathcal I}}
\def\cJ{{\mathcal J}}
\def\cK{{\mathcal K}}
\def\cL{{\mathcal L}}
\def\cM{{\mathcal M}}
\def\cN{{\mathcal N}}
\def\cO{{\mathcal O}}
\def\cP{{\mathcal P}}
\def\cQ{{\mathcal Q}}
\def\cR{{\mathcal R}}
\def\cS{{\mathcal S}}
\def\cT{{\mathcal T}}
\def\cU{{\mathcal U}}
\def\cV{{\mathcal V}}
\def\cW{{\mathcal W}}
\def\cX{{\mathcal X}}
\def\cY{{\mathcal Y}}
\def\cZ{{\mathcal Z}}

\newcommand{\comm}[1]{\marginpar{%
\vskip-\baselineskip 
\raggedright\footnotesize
\itshape\hrule\smallskip#1\par\smallskip\hrule}}





\def\ve{\varepsilon}

\hyphenation{re-pub-lished}

\def\ord{{\mathrm{ord}}}
\def\Nm{{\mathrm{Nm}}}
\renewcommand{\vec}[1]{\mathbf{#1}}

\def \F{{\bbbf}}
\def \L{{\bbbl}}
\def \K{{\bbbk}}
\def \Z{{\bbbz}}
\def \N{{\bbbn}}
\def \Q{{\bbbq}}
\def\E{{\mathbf E}}
\def\bH{{\mathbf H}}
\def\G{{\mathcal G}}
\def\O{{\mathcal O}}
\def\cS{{\mathcal S}}
\def \R{{\bbbr}}
\def\Fp{\F_p}
\def \fp{\Fp^*}
\def\\{\cr}
\def\({\left(}
\def\){\right)}
\def\fl#1{\left\lfloor#1\right\rfloor}
\def\rf#1{\left\lceil#1\right\rceil}

\def\Zm{\Z_m}
\def\Zt{\Z_t}
\def\Zp{\Z_p}
\def\Um{\cU_m}
\def\Ut{\cU_t}
\def\Up{\cU_p}

\def\ep{{\mathbf{e}}_p}

\def\e{{\mathbf{e}}}

\def \Prob{{\mathrm {}}}

\def\LC{{\cL}_{C,\cF}(Q)}
\def\LCn{{\cL}_{C,\cF}(nG)}
\def\Mrs{\cM_{r,s}\(\F_p\)}

\def\Fbar{\overline{\F}_q}
\def\Fn{\F_{q^n}}
\def\En{\E(\Fn)}

\def \hatf{\widehat{f}}

\def\mand{\qquad \mbox{and} \qquad}

\def\MOV{{\bf{MOV}}}

\title{\bf On pseudosquares and pseudopowers}

\author{
{\sc Carl Pomerance}\\
{Department of Mathematics}\\ {Dartmouth College}\\
{Hanover, NH 03755-3551, USA} \\
{\tt carl.pomerance@dartmouth.edu} \\
\and
{\sc Igor E. Shparlinski}\\
                 {Department of Computing}\\
{Macquarie University}\\
{ Sydney, NSW 2109,
Australia}\\
{\tt igor@ics.mq.edu.au}
}

\date{}

\maketitle

\begin{abstract} 
Introduced by Kraitchik and Lehmer,
an $x$-pseudosquare is a positive integer $n\equiv1\pmod 8$ that
is a quadratic residue for each odd prime $p\le x$, yet is not a square.
We use bounds of character sums 
to prove that pseudosquares are equidistributed
in fairly short intervals. 
An $x$-pseudopower to base $g$ is a positive integer which is not a power
of $g$ yet is so modulo $p$ for all primes $p\le x$.
It is conjectured by
Bach, Lukes, Shallit, and Williams
that the least such number is at most $\exp(a_g x/\log x)$
for a suitable constant $a_g$.
A bound of $\exp(a_g x\log\log x/\log x)$ is proved
conditionally on the Riemann Hypothesis for Dedekind
zeta functions, thus improving on a recent conditional exponential
bound of Konyagin and the present authors.  We also give
a GRH-conditional equidistribution result for pseudopowers
that is analogous to our unconditional result for pseudosquares.
\end{abstract}


\section{Introduction}

\subsection{Pseudosquares}
\label{ssec1.1}

An $x$-pseudosquare is a nonsquare positive integer
$n$ such that $n\equiv1\pmod8$ and $(n/p)=1$ for each
odd prime $p\le x$.
The subject of pseudosquares was initiated by Kraitchik
and more formally by Lehmer in~\cite{L}.  It was later
shown by Weinberger (see~\cite{W}) that
if the Generalized Riemann Hypothesis (GRH) holds, then
the least $x$-pseudosquare, call it $N_x$, satisfies
$N_x\ge \exp(cx^{1/2})$ for a positive constant $c$.
The interest in this inequality is that there is a
primality test, due to Selfridge and Weinberger,
for integers $n<N_x$ that requires
the verification of some simple Fermat-type congruences
for prime bases
$p\le x$.  Thus, a large lower bound
for $N_x$ leads to a fast primality test, and in particular
this result gives an alternate and somewhat simpler form
of Miller's GRH-conditional polynomial-time deterministic
primality test.  See~\cite{W} for details.

By the GRH, we mean the Riemann Hypothesis for Dedekind zeta
functions, that is, for algebraic number fields.  Note that
this conjecture subsumes the Extended Riemann Hypothesis (ERH),
which is the Riemann Hypothesis for rational Dirichlet $L$-functions.
The Weinberger lower bound for $N_x$ in fact only requires the ERH.

As the concept of an $x$-pseudosquare is a natural one, it is also of
interest to find a reasonable upper bound for $N_x$ and also to
study the distribution of $x$-pseudosquares.  Let $M(x)$ denote the
product of the primes up to $x$.  For nonsquare integers $n$ coprime to $M(x)$,
the ``probability'' that $n$ satisfy $n\equiv1\pmod 8$ and $(n/p)=1$
for all odd primes $p\le x$ is $2^{-\pi(x)-1}$. Thus, it is reasonable
perhaps to conjecture that
$$
N_x\le 2^{(1+o(1))\pi(x)},
$$
for example, see Bach and Huelsbergen~\cite{BH}.  
In~\cite{S3}, Schinzel proves conditionally on the GRH that
\begin{equation}
\label{eq:Conj}
N_x\le 4^{(1+o(1))\pi(x)},
\end{equation}
and in particular, he conditionally shows that this inequality holds 
as well for the smallest
prime $x$-pseudosquare.  Unconditionally, he uses
the Burgess bound~\cite{B} (see also~\cite[Theorem~12.6]{IK})
to show that
\begin{equation}
\label{eq:Schinzel}
N_x\le\exp((1/4+o(1))x).
\end{equation}

We start with an observation, communicated to us by K.~Soundararajan,
that the pigeonhole
principle (used in the same fashion as in~\cite[Lemma~10.1]{GrSo}) 
gives an unconditional proof of~\eqref{eq:Conj}, though not for
prime pseudosquares.  Indeed, let us put $X =2^{\pi(x)} x$ and consider the
$$
(1+ o(1)\frac{X}{ \log X} =  (1+ o(1)\frac{2^{\pi(x)} \log x}{  \log 2}
$$
vectors of Legendre symbols
$$
\(\(\frac{\ell}{3}\), \(\frac{\ell}{5}\), \ldots , \(\frac{\ell}{p}\)\)
$$
for all primes $\ell \in (x, X]$,
where $p$ is the largest prime with $p \le x$. Clearly
there are at most $2^{\pi(x)-1}$ possbilities for such vectors,
so for large $x$ there are five distinct primes 
$\ell_1,\dots,\ell_5 \in (x, X]$
for which they coincide. Thus, at least two of them, say $\ell_1,\ell_2$,
have the property that $\ell_1\equiv\ell_2\pmod 8$.
Then $\ell_1\ell_2$ is an $x$-pseudosquare
and we have $N_x \le \ell_1\ell_2 \le X^2$, implying~\eqref{eq:Conj}.
(Note that it is not necessary that the numbers $\ell$ be prime
in this proof, just coprime to $M(x)$.)

Our contribution to the subject of pseudosquares is on their 
equidistribution.  
For this we follow Schinzel's proof of~\eqref{eq:Schinzel}, but
use a character sum estimate given in~\cite[Corollary 12.14]{IK},
which dates back to work of Graham and Ringrose~\cite{GR},
to prove the following result.
Let $\cS_x$ be the set of $x$-pseudosquares.

\begin{theorem}
\label{squaretheorem}
For any interval $(A, A+N] \subseteq (0,\infty)$, uniformly over
$N\ge\exp(3x/\log\log x)$,   we have
\[
\# \(\cS_x \cap (A, A+N]\) = (1+o(1))  \frac{N}{M(x)} \#   \(\cS_x
\cap (0, M(x)]\), \qquad
x \to \infty.
\]
\end{theorem}

We also show
\begin{equation}
\label{eq:Card Qx}
\#   \(\cS_x \cap (0, M(x)]\) = (1+o(1))
\frac{M(x)}{2^{\pi(x)+1}e^\gamma\log x}
, \qquad x \to \infty,
\end{equation}
so that
one can rewrite Theorem~\ref{squaretheorem} in a more explicit form.

Note that Granville and Soundararajan \cite{GrSo} also discuss the
equidistibution of $x$-pseudosquares via the Graham--Ringrose
result on character sums, but their context is different
and it is not clear that Theorem \ref{squaretheorem} follows
directly from their paper.

\subsection{Pseudopowers}

Let $g$ be a fixed integer with $|g|\ge 2$.
Following Bach, Lukes, Shallit, and
Williams~\cite{BLSW}, we say that  an
integer $n>0$ is  an {\it $x$-pseudopower to  base $g$\/}
if $n$ is not a power of
$g$ over the integers but is
a power of $g$  modulo all primes $p\le x$.
Denote by $q_g(x)$ the least $x$-pseudopower to base $g$.

In~\cite{BLSW} it is conjectured that for each fixed $g$,
there is a number $a_g$ such that for $x\ge 2$,
\begin{equation}
\label{BLSWconj}
q_g(x)\le\exp(a_gx/\log x).
\end{equation}
In addition, a heuristic argument is given for~\eqref{BLSWconj},
with numerical evidence presented
in the case of $g=2$.
For any $g$, we have (see~\cite{KPS})
the trivial bound $q_g(x)\le 2M(x)+1$, where $M(x)$
is the product of the primes $p\le x$.  Thus,
\[
q_g(x)\le \exp((1+o(1))x).
\]
Using an estimate for exponential sums
due to Heath-Brown and Konyagin~\cite{HBK} and results
of Baker and Harman~\cite{BH1,BH2} on the Brun--Titchmarsh theorem
on average, Konyagin, Pomerance, and Shparlinski~\cite{KPS}
proved that
\[
q_g(x)\le\exp(0.88715x)
\]
for all sufficiently large $x$ and all integers $g$ with
$2\le|g|\le x$.  Further,
it was noted in~\cite{KPS}
that the method implied that for fixed $g$,
\[
q_g(x)\le \exp((1/2+o(1))x),
\]
assuming the GRH.  In this paper we make further
progress towards~\eqref{BLSWconj}, again assuming the GRH.
Our proof makes use of the approach in
Schinzel~\cite{S3} for pseudosquares.

\begin{theorem}
\label{qgx bound}
Assume the GRH.
Then for each fixed integer $g$ with $|g|\ge 2$ there is
a number $a_g$ such that for $x\ge 3$,
\[
q_g(x)\le\exp(a_g x\log\log x/\log x).
\]
\end{theorem}

We are also able to prove an equidistribution result conditional on the GRH
that is similar in strength to Theorem~\ref{squaretheorem}.
Let $\cP_x$ be the set of $x$-pseudopowers base $g$.

\begin{theorem}
\label{equipower}
Assume the GRH.  Let $g$ be a fixed integer with $|g|>1$.  There is
a positive number $b_g$ such that
for any interval $(A, A+N]\subseteq (0,\infty)$, uniformly over
$N\ge\exp(b_gx/\log\log x)$,   we have
\[
\# \(\cP_x \cap(A, A+N]\) = (1+o(1))
\frac{N}{M(x)} \# \(\cP_x\cap(0, M(x)]\), \qquad
x \to \infty.
\]
\end{theorem}

We derive an asymptotic formula for  $\# \(\cP_x \cap(0, M(x)]\)$
in Section~\ref{sec-equipower}, see~\eqref{eq:Card Wx},
so that one can get a more explicit form of Theorem~\ref{equipower}.

We note that in~\cite{KPS} an unconditional version
of Theorem~\ref{equipower}  is given which however
holds only for $N \ge \exp(0.88715x)$. Under the GRH, the method
of~\cite{KPS} gives a
somewhat stronger result but still
requires $N$ to be rather large, namely it applies only
to $N \ge \exp\((0.5 + \varepsilon) x\)$ for an arbitrary
$\varepsilon>0$.

As for lower bounds for $q_g(x)$, it follows from
Schinzel~\cite{S1,S2} that
\[
q_g(x)\to\infty, \qquad x\to\infty.
\]
In~\cite{BLSW}
it is shown that assuming the GRH
there is a number $c_g>0$ such that
\[
q_g(x)\ge\exp(c_g\sqrt x(\log\log x)^3/(\log x)^2).
\]

\subsection{Notation}

We recall that the notation
$U = O(V)$ and $U\ll V$ are equivalent
to the assertion that the inequality $|U|\le c\,V$ holds for some
constant $c>0$.

\section{Distribution of pseudosquares}

In this section we   prove
Theorem~\ref{squaretheorem}
by making use of the following character sum estimate,
which is~\cite[Corollary 12.14]{IK}.

\begin{lemma}
\label{IKlemma}
Let $\chi$ be a primitive character to the squarefree modulus $q>1$.
Suppose all prime factors of $q$ are at most $N^{1/9}$ and let
$r$ be an integer with $N^r\ge q^3$.  Then for any number $A$,
\[
\left|\sum_{A<n\le A+N}\chi(n)\right|
\le 4N\tau(q)^{r/2^r}q^{-1/r2^r},
\]
where $\tau(q)$ is the number of positive divisors of $q$.
\end{lemma}

Recall that $M(x)$ is the product of the primes in $[1,x]$.
Let $x$ be a large number and
let $\overline \cS_x$ denote the set of positive integers
$n\equiv1\pmod 8$ with $(n/p)=1$ for each
odd prime $p\le x$.
That is, $\overline \cS_x$ consists of the $x$-pseudosquares and actual
squares coprime to $M(x)$.
In particular, $\cS_x\subseteq \overline \cS_x$.
We let $M_2(x)=M(x)/2$,
the product of the odd primes up to $x$.

Theorem~\ref{squaretheorem} is routine once $N$ is large compared
with $M(x)$, so we  assume that $N\le M(x)^2$.
Note that for
a positive integer
$n$ with $(8n+1,M_2(x))=1$, we have that
\begin{equation}
\label{squarechar}
\prod_{p\mid M_2(x)}\(1+\(\frac{8n+1}{p}\)\) =
\left\{\begin{array}{ll}
2^{\pi(x)-1}, &\quad\text{if $8n+1\in\overline \cS_x$;}\\
0,&\quad\text{else.}
\end{array}\right.
\end{equation}
Thus, if $A,N$ are positive numbers, then the sum
\[
S_{A,N}:=
\sum_{\substack{A<8n+1\le A+N\\(8n+1,M_2(x))=1}}
\prod_{p\mid M_2(x)}\(1+\(\frac{8n+1}{p}\)\)
\]
satisfies
\begin{equation}
\label{squarecount}
S_{A,N}=2^{\pi(x)-1}\#(\overline \cS_x\cap(A,A+N]).
\end{equation}

The product in~\eqref{squarechar} can be
expanded, so that we have
\begin{equation}
\label{quadsum}
\begin{split}
S_{A,N} & =~\sum_{\substack{A<8n+1\le A+N\\(8n+1,M_2(x))=1}}
\sum_{f\mid M_2(x)}\(\frac{8n+1}{f}\)\\
&=~\sum_{f\mid M_2(x)}\sum_{\substack{A<8n+1\le A+N\\ (8n+1,M_2(x))=1}}
\(\frac{8n+1}{f}\).
\end{split}
\end{equation}
The contribution to $S_{A,N}$ from $f=1$ is
\begin{equation}
\label{mainterm}
\sum_{\substack{A<8n+1\le A+N\\ (8n+1,M_2(x))=1}}1 \sim
\frac{N}{4e^\gamma\log x}
\end{equation}
uniformly for $A,N$ with $N\ge \exp(x^{1/2})$.  This estimate follows
immediately from the fundamental lemma of the sieve; for example,
see~\cite[Theorem 2.5]{HR}.

Suppose now that $f\mid M_2(x)$, $f>1$ is fixed.  We can rewrite the
contribution in~\eqref{quadsum} corresponding to $f$ as
\begin{equation}
\label{squarecharsum2}
R_f= \sum_{d\mid M_2(x)/f}\mu(d)\sum_{\substack{A<8n+1\le A+N\\ d\mid 8n+1}}
\(\frac{8n+1}{f}\),
\end{equation}
where $\mu(d)$ is the M{\"o}bius function.

The P\'olya--Vinogradov inequality (see~\cite[Theorem~12.5]{IK})
immediately implies that
\begin{equation}
\label{pv}
|R_f|  \le 3\cdot2^{\pi(x)-2}\sqrt{f}\log f<2^{\pi(x)}\sqrt{f}\log f
\end{equation}
for any choice of $f>1$. We use~\eqref{pv} when $f$ is
not much larger than $N$,
namely we use it
when
\[
f\le N^{r2^r/(r2^{r-1}+1)},
\]
where $r$ shall be chosen later.  In this case it gives
\begin{equation}
\label{eq:Rf bound}
|R_f| \le  2^{\pi(x)}N^{1-2/(r2^r+2)}\log(N^2)
\le 2^{(1 + o(1))\pi(x)}N^{1-2/(r2^r+2)} .
\end{equation}

For large values of $f$, that is,
when
\begin{equation}
\label{eq:break}
f>N^{r2^r/(r2^{r-1}+1)},
\end{equation}
     we use a different approach which relies
on Lemma~\ref{IKlemma}.

Let $r_f=(1-f^2)/8$, so that $r_f$ is an integer
and $8r_f\equiv 1\pmod f$.
Then
\begin{eqnarray*}
R_f & = &
\sum_{d\mid M_2(x)/f}\mu(d)\(\frac df\)
\sum_{\substack{A<dk\le A+N\\ k\equiv d\kern-5pt\pmod{8}}}\(\frac kf\)\\
& = &
    \sum_{d\mid M_2(x)/f}\mu(d)\(\frac{8d}f\)
\sum_{A<d(8l+d)\le A+N}
\(\frac{l+dr_f}{f}\)\\
& = &
\sum_{d\mid M_2(x)/f}\mu(d)\(\frac{8d}f\)
\sum_{m\in\cI_{d,f}}
\(\frac{m}{f}\),
\end{eqnarray*}
where $\cI_{d,f} =[B_{d,f}+1, B_{d,f}+N_{d,f}]$,
an interval of length
\[
N_{d,f} = \frac{N}{8d} + O(1).
\]
Thus,
\begin{equation}
\label{eq:Rf}
|R_f|  \le
     \sum_{d\mid M_2(x)/f}\left| \sum_{m \in \cI_{d,f}}
\(\frac{m}{f}\)  \right|.
\end{equation}

The character sums  in~\eqref{eq:Rf} where $8d>N^{0.1}$ are trivially bounded
by $N_{d,f} = O(N^{0.9})$ in absolute value, so their total
contribution to $R_f$ is
\begin{equation}
\label{larged}
\sum_{\substack{d\mid M_2(x)/f\\ 8d > N^{0.1}}}\left| \sum_{m \in \cI_{d,f}}
\(\frac{m}{f}\)  \right|  \ll 2^{\pi(x)}N^{0.9}.
\end{equation}

We now assume that $8d\le N^{0.1}$.
Note that the conductors $f$ of the characters
which appear in~\eqref{eq:Rf} are squarefree.
We  choose $r$ as the largest integer with
\[
r2^r+2\le\frac{ \log x}{\log\log x}
\]
and apply Lemma~\ref{IKlemma} to the inner sum
in~\eqref{eq:Rf} with this value of $r$.  To do this we
need
\[
(N/(8d))^r\ge f^3 \mand x\le (N/(8d))^{1/9}.
\]
These inequalities hold   since
\begin{equation}
\label{eq:r and N}
r= \(\frac{1}{\log 2}+   o(1)\)\log\log x \mand
N \ge \exp(3x/\log\log x),
\end{equation}
so that
\begin{equation}
\label{eq:Cond 1}
\(\frac{N}{8d}\)^r\ge N^{0.9r}\ge \exp(2.7rx/\log\log x)
\ge M(x)^3 \ge  f^3
\end{equation}
and
\begin{equation}
\label{eq:Cond 2}
\(\frac{N}{8d}\)^{1/9}\ge N^{0.1}\ge \exp(0.3x/\log\log x) \ge x.
\end{equation}
for all large $x$.

Thus, by Lemma~\ref{IKlemma},
\begin{equation}
\label{smalld}
\sum_{\substack{d\mid M_2(x)/f\\ 8d \le  N^{0.1}}}\left| \sum_{m \in \cI_{d,f}}
\(\frac{m}{f}\)  \right| \le 4\cdot2^{\pi(x)-1}N2^{(\pi(x)-1)r/2^r}f^{-1/r2^r}.
\end{equation}
We now derive
from~\eqref{larged},~\eqref{smalld},~and \eqref{eq:break} that
\[
|R_f| \le 2^{(1+o(1))\pi(x)}N^{0.9} +  2^{(1 + o(1))\pi(x)}N f^{-1/r2^r}
\le 2^{(1+o(1))\pi(x)}N^{1-2/(r2^r+2)}
\]
which is also the bound in~\eqref{eq:Rf bound}
for those $f>1$ not satisfying~\eqref{eq:break}.

Now summing on $f$ we see that the contribution to
$S_{A,N}$ from values of $f>1$ is at most
\begin{equation*}
\begin{split}
\sum_{\substack{f \mid M_2(x)\\ f>1}} |R_f |&\le~
4^{(1+o(1))\pi(x)}N^{1-2/(r2^r+2)}\\
&=~N^{1-2/(r2^r+2)}\exp((\log
4+o(1))x/\log x).
\end{split}
\end{equation*}
Note that by our choice of $r$ we have
\[
N^{2/(r2^r+2)}\ge N^{2  \log\log x/\log x}\ge\exp(6  x/\log x).
\]
Since $6  >\log 4$, we have that the contribution to~\eqref{quadsum}
from terms with $f>1$ is small compared to the main term given
by~\eqref{mainterm}, so that
\[
S_{A,N} = (1 + o(1))\frac{N}{4e^\gamma\log x}.
\]
Together with~\eqref{squarecount}, we now have
\[
\#(\overline \cS_x\cap(A,A+N]) = (1 + o(1))
\frac{N}{2^{\pi(x)+1}e^\gamma\log x}.
\]
Since the number of squares in the interval $(A,A+N]$ is at
most $N^{1/2}$,
we obtain
\[
\# \(\cS_x \cap (A, A+N]\) = (1 + o(1)) \frac{N}{2^{\pi(x)+1}e^\gamma\log x}.
\]
Taking $A =0$ and $N = M(x)$ we derive~\eqref{eq:Card Qx}
and conclude the proof of Theorem~\ref{squaretheorem}.

\medskip
We remark that by being a little more careful with the estimates,
we can prove the theorem with ``3" replaced with any fixed number
larger than $\log 8$.

\section{Distribution of pseudopowers}

\subsection{Proof of Theorem~\ref{qgx bound}}
\label{sec-smallpower}

Let $g$ be a given integer with $|g|\ge 2$ which we assume
to be fixed.  Let $p_g(x)$ be the least positive integer
which is not a power of $g$ yet is a power of $g$ modulo
every prime $p\le x$ with $p\nmid g$.  It is easy to see
that
\begin{equation}
\label{easytosee}
q_g(x)\le gp_g(x).
\end{equation}
Indeed, the integer $gp_g(x)$ is not a power of $g$, it
is a power of $g$ modulo every prime $p\le x$ with $p\nmid g$,
and it is zero modulo $p$ for every prime $p\mid g$, and so is a power
of $g$ modulo these primes too.

For every prime $p\nmid g$ let $l_g(p)$ be the multiplicative order of $g$
modulo $p$, and let $i_g(p)=(p-1)/l_g(p)$, the index of the subgroup
of powers of $g$ in the multiplicative group modulo $p$.  Let
\[
M_g(x)=\prod_{\substack{p\le x\\ p\nmid g}}p,\quad
I_g(x)=\prod_{p\mid M_g(x)}i_g(p).
\]
It follows from~\cite[Theorem~1]{KPS} that
$I_g(x)\le\exp(0.42x)$ for all sufficiently large $x$.  We conditionally
improve this result.

\begin{lemma}
\label{Igx bound}
Assume the GRH.
There is a number $c_g$ such that for $x\ge3$,
\[
I_g(x)\le\exp(c_g x\log\log x/\log x)
\]
\end{lemma}

\begin{proof}
In Kurlberg and Pomerance~\cite[Theorem 23]{KP},
following ideas of Hooley~\cite{H} and Pappalardi~\cite{P},
it is shown conditionally on the GRH that
\[
\sum_{\substack{p\mid M_g(x)\\ l_g(p)\le p/y}}1
\ll_g \frac{\pi(x)}{y}+\frac{x\log\log x}{(\log x)^2}
\]
for $1\le y\le\log x$.  Applying this result with $y=\log x$,
we have
\[
\sum_{\substack{p\mid M_g(x)\\ i_g(p)\ge \log x}}1
\ll_g\frac{x\log\log x}{(\log x)^2}.
\]
Indeed, $i_g(p)\ge y$ implies that $l_g(p)\le(p-1)/y<p/y$.
Since we trivially have $i_g(p)\le x$ for each prime $p\mid M_g(x)$
with $p\nmid g$, we thus have
\[
I_g(x)=\prod_{\substack{p\mid M_g(x)\\ i_g(p)<\log x}}i_g(p)
\prod_{\substack{p\mid M_g(x)\\ i_g(p)\ge\log x}}i_g(p)
\le(\log x)^{\pi(x)}x^{O_g(x\log\log x/(\log x)^2)}.
\]
The lemma follows.
\end{proof}

For each prime $p\nmid g$, let $\chi_p$ be a character modulo $p$
of order $i_g(p)$.  Then
\[
\sum_{j=1}^{i_g(p)}\chi_p^j(n)
=
\left\{ \begin{array}{ll}
i_g(p), & \quad \text{if}\ n\ \text{is a power of}\ g\ \text{modulo}\ p,\\
0, & \quad \text{else}.
\end{array} \right.
    \]
Thus,
\begin{equation}
\label{charprod}
\begin{split}
\prod_{p\mid M_g(x)}
\sum_{j=1}^{i_g(p)}&\chi_p^j(n)\\
= &
\left\{ \begin{array}{ll}
I_g(x), & \quad \text{if $n$ is a power of $g$
modulo every $p\mid M_g(x)$},\\ 0, & \quad \text{else}.
\end{array} \right.
\end{split}
\end{equation}
Let $\Lambda(n)$ denote the von Mangoldt function.
    From the definition of $p_g(x)$ we
deduce that
\[
S_g:= \sum_{n<p_g(x)}\Lambda(n)
\prod_{p\mid M_g(x)}
\sum_{j=1}^{i_g(p)}\chi_p^j(n)
=I_g(x)\sum_{\substack{n<p_g(x)\\ n\,{\rm is\,a\,power\,of}\,g}}\Lambda(n).
\]
The last sum is 0 if $g$ is not a prime or prime power, and in
any event is always at most $\log p_g(x)\ll x$.  Thus,
\begin{equation}
\label{Sbound}
S_g\ll I_g(x)x.
\end{equation}

We now multiply out the product in~\eqref{charprod}; it is
seen as a sum of $I_g(x)$ characters modulo $M_g(x)$.  The
contribution to $S_g$ from the principal character
$\prod_p\chi_p^{i_g(p)}$ is $\psi(p_g(x))+O(x)$.  We may
assume that $p_g(x)\ge e^{\pi(x)}$ since otherwise the
theorem follows immediately from~\eqref{easytosee}, so that
the contribution to $S_g$ from the principal character is
$(1 + o(1))p_g(x)$, by the prime number theorem.

The contribution to $S_g$ from each nonprincipal character
$\chi$ is
\[
\sum_{n<p_g(x)}\chi(n)\Lambda(n)
\ll p_g(x)^{1/2}\((\log M_g(x))^2+(\log p_g(x))^2\)\ll p_g(x)^{1/2}x^2
\]
assuming the GRH.  Hence, the contribution to $S_g$ from
nonprincipal characters is
$O(I_g(x)p_g(x)^{1/2}x^2)$.
Thus,
\[
S_g=(1+o(1))p_g(x)+O(I_g(x)p_g(x)^{1/2}x^2),
\]
so that from~\eqref{Sbound} we deduce that
\[
p_g(x)\ll I_g(x)^2x^4.
\]
Theorem~\ref{qgx bound} now follows immediately from~\eqref{easytosee}
and Lemma~\ref{Igx bound}.

\medskip

We remark that an alternate way to handle primes dividing $g$
is to eschew~\eqref{easytosee} and instead multiply the product
in~\eqref{charprod} by $\sum_{\chi \bmod g}\chi(n)$.  This sum
is $\varphi(g)$ when $n\equiv1\pmod g$ and is 0 otherwise.
Note that a number that is 1~mod~$p$ is always a power of $g$ modulo $p$.
Although this is somewhat more complicated, it does lead to a proof
that there is a {\it prime} number below the bound
$\exp(a_gx\log\log x/\log x)$ that is
an $x$-pseudopower base $g$.

\subsection{Proof of Theorem~\ref{equipower}}
\label{sec-equipower}

We use
the method of proof of Theorem~\ref{squaretheorem}.
Accordingly we only outline some new elements
and   suppress the details.

For a nonzero integer $n$,
let $\rad(n)$, the {\it radical\/}
of $n$, be the largest squarefree divisor of $n$.  That is, $\rad(n)$
is the product of the distinct prime factors of $n$.  Also,
let
$\omega(n)$ denote the number of distinct prime factors of $n$.
Let $\overline \cP_x$ denote the set of positive integers which are
either an $x$-pseudopower base $g$ or a true power of $g$.
Then a positive integer $n\in\overline \cP_x$ if and only if both
\begin{enumerate}
\item[(i)] $n$ is in the subgroup
$\langle g\rangle$ of $(\Z/p\Z)^*$ when $p\le x$ and $p\nmid g$;
\item[(ii)] $n\equiv0$ or $1\pmod p$ when $p\le x$ and $p\mid g$.
\end{enumerate}
Assuming then that $x\ge |g|$,
the cardinality of $\overline \cP_x\cap(0,M(x)]$ is
\begin{equation*}
\begin{split}
2^{\omega(g)}\prod_{p\mid M_g(x)}l_g(p)
=2^{\omega(g)}\prod_{p\mid M_g(x)}\frac{p-1}{i_g(p)}
=&\frac{2^{\omega(g)}\varphi(M_g(x))}{I_g(x)}\\
\sim&\frac{2^{\omega(g)} M(x)}{e^\gamma\varphi(\rad(g))I_g(x)\log x},
\end{split}
\end{equation*}
by the formula of Mertens.

It is easy to see that this expression is
exponentially large,  either from the observation
that $l_g(p)\ge2$ whenever $\(\frac gp\)=-1$, so that
$\#(\overline \cP_x\cap(0,M(x)]))\ge 2^{(1/2+o(1))\pi(x)}$, or using
$I_g(x)\le e^{0.42x}$ from~\cite{KPS}.  Further, the number of true
powers of $g$ in $(0,M(x)]$ is small; it is $O(x)$.  Thus,
\begin{equation}
\label{eq:Card Wx}
\# \(\cP_x \cap(0, M(x)]\) = (1+ o(1))
\frac{2^{\omega(g)} M(x)}{e^\gamma\varphi(\rad(g))I_g(x)\log x}
, \qquad x \to \infty.
\end{equation}

To prove Theorem~\ref{equipower} we
again use Lemma~\ref{Igx bound}, which is GRH-conditional.
But the framework of the proof follows the argument of
Theorem~\ref{squaretheorem}, and in particular it uses the unconditional
Lemma~\ref{IKlemma}.
Notice that the proof of Theorem~\ref{qgx bound} used Riemann Hypotheses
a second time, namely in the estimation of the weighted character sums.
Now we use unweighted character sums and
so are able to use Lemma~\ref{IKlemma}.  The set-up is as follows.
Let
\begin{equation}
\label{equichar}
P_{A,N} = \sum_{ab=\rad(g)}
\sum_{\substack{A<n\le A+N\\
n\equiv0\kern-7pt\pmod{a}\\n\equiv1\kern-7pt\pmod{b}}}
\prod_{p\mid M_g(x)}\sum_{j=1}^{i_g(p)}\chi_p^j(n).
\end{equation}
This expression counts integers $n\in(A,A+N]$ that
are 0 or $1\pmod p$ for each prime $p\mid g$ and in the subgroup
$\langle g\rangle$
of $(\Z/p\Z )^*$ for each prime $p\le x$ with $p\nmid g$.
Namely, it counts members of $\overline \cP_x$, and does so with the
weight $I_g(x)$.

To prove Theorem~\ref{equipower} one then expands
the product in~\eqref{equichar}.  As usual, the contribution from the
principal character
is easily estimated: it is the number of integers $n\in(A,A+N]$
which are 0 or $1\pmod p$ for each prime $p\mid g$ and are coprime
to $M_g(x)$.  Thus, the principal character gives the contribution
\[
(1+o(1)) \frac{2^{\omega(g)}}{\rad(g)}\cdot
\frac{\varphi(M_g(x))}{M_g(x)}N
=(1+o(1)) \frac{2^{\omega(g)}N}{e^\gamma\varphi(\rad(g))\log x},
\]
which when divided by the weight $I_g(x)$ gives the main term for our count.

The nonprincipal characters have conductors corresponding to
those divisors $f$ of $M_g(x)$ with $f>1$ and $i_g(p)>1$ for
each prime $p\mid f$.  For such integers $f$, the
characters that occur with conductor $f$ are induced by characters in the
set
\[
X_f=\Bigl\{\prod_{p\mid f}\chi_p^{j_p}~:~1\le j_p\le i_g(p)-1\hbox{
for }p\mid f\Bigr\}.
\]
Thus, the contribution of the nonprincipal characters to $P_{A,N}$ is
\[
P^*_{A,N}:=\sum_{ab=\rad(g)}\sum_{\substack{f\mid M_g(x)\\ f>1}}
\sum_{\chi\in X_f}
\sum_{\substack{A<n\le A+N\\ n\equiv0\kern-5pt\pmod{a}\\
n\equiv1\kern-5pt\pmod{b}\\(n,M_g(x)/f)=1}}\chi(n),
\]
where $X_f$ is empty if $i_g(p)=1$ for some prime $p\mid f$.
The inner sum is
\[
\sum_{d\mid M_g(x)/f}\mu(d)
\sum_{\substack{A<n\le A+N\\ n\equiv0\kern-5pt\pmod{ad}\\
n\equiv1\kern-5pt\pmod{b}}}\chi(n).
\]
As before we estimate the character sum here trivially if $d$ is large,
we use the P\'olya--Vinogradov inequality if $f$ is small,
and we use Lemma~\ref{IKlemma} in the remaining cases.
However, we modify slightly the choice of $r$
in the proof of Theorem~\ref{squaretheorem} and take it now
as the largest integer with
\[
r2^r+2\le\frac{ \log x}{(\log\log x)^2}.
\]
Then we still have~\eqref{eq:r and N} and thus the conditions~\eqref{eq:Cond 1}
and~\eqref{eq:Cond 2} are still
satisfied so that Lemma~\ref{IKlemma} may be used.
Since each $|X_f|<I_g(x)$,
we obtain
\[
|P^*_{A,N}|\le 4^{(1+o(1))\pi(x)}I_g(x)N^{1-2/(r2^r+2)} .
\]
This leads us to the asymptotic formula
\[
P_{A,N}=
(1+o(1)) \frac{2^{\omega(g)}N}{e^\gamma\varphi(\rad(g))\log x}
+ O\(I_g(x)4^{(1+o(1))\pi(x)}N^{1-2/(r2^r+2)} \) .
\]

For $N\ge\exp(b_gx/\log\log x)$
we have
\[
N^{2/(r2^r+2)}\ge N^{2   (\log\log x)^2/\log x}
>\exp\(2   b_g x\log \log x/\log x\).
\]
Using Lemma~\ref{Igx bound} and taking $b_g =   c_g$, we
conclude  the
proof of Theorem~\ref{equipower}.

\section*{Acknowledgments}
The authors are grateful to K.~Soundararajan for communicating the proof
we presented in Section~\ref{ssec1.1} and for his permission to give it
here.

The first author was supported in part by NSF grant DMS-0703850.
The second author was supported in part by ARC grant DP0556431.

\end{document}